\documentclass[onecolumn]{elsart3p}
\usepackage{amssymb}

\usepackage[english,francais]{babel}

\newtheorem{theorem}{Theorem}[section]

\newtheorem{e-proposition}[theorem]{Proposition}

\newtheorem{e-definition}[theorem]{Definition\rm}

\newtheorem{theoreme}{Th\'eor\`eme}

\newcommand{\R}{\mathbb{R}}

\newcommand{\href}[1]{(\ref{#1})}

\newcommand{\eqref}[1]{(\ref{#1})}
\newcommand{\be}[1]{\begin{equation}\label{#1}}
\newcommand{\ee}{\end{equation}}
\newenvironment{multline*}{\begin{eqnarray*}}{\end{eqnarray*}}
\newenvironment{equation*}{\[}{\]}

\def\og{\leavevmode\raise.3ex\hbox{$\scriptscriptstyle\langle\!\langle$~}}
\def\fg{\leavevmode\raise.3ex\hbox{~$\!\scriptscriptstyle\,\rangle\!\rangle$}}

\journal{the Acad\'emie des sciences}
\begin{document}

\centerline{Partial Differential Equations/Probability}
\begin{frontmatter}

\selectlanguage{english}


\title{About Kac's Program in Kinetic Theory}

\selectlanguage{english}
\author[Ceremade]{St\'ephane Mischler},
\ead{mischler@ceremade.dauphine.fr}
\author[Cambridge]{Cl\'ement Mouhot},
\ead{C.Mouhot@dpmms.cam.ac.uk}

\address[Ceremade]{Ceremade (UMR CNRS no. 7534), Universit\'e
  Paris-Dauphine, Place de-Lattre-de-Tassigny, 75775 Paris cedex 16,
  France} \address[Cambridge]{DPMMS, Centre for Mathematical Sciences,
  University of Cambridge, Wilberforce Road, Cambridge CB2 0WA, UK}

\medskip
\begin{center}
{\small Received \dots 2011; \dots \\ Presented by \dots}
\end{center}

\begin{abstract}
  \selectlanguage{english} 

  In this Note we present the main results from the recent
  work~\cite{mm-article}, which answers several conjectures raised
  fifty years ago by Kac \cite{kac}. There Kac introduced a
  many-particle stochastic process (now denoted as Kac's master
  equation) which, for chaotic data, converges to the spatially
  homogeneous Boltzmann equation. We answer the three following
  questions raised in \cite{kac}: (1) prove the propagation of chaos
  for realistic microscopic interactions (i.e. in our results: hard
  spheres and true Maxwell molecules); (2) relate the time scales of
  relaxation of the stochastic process and of the limit equation by
  obtaining rates independent of the number of particles; (3) prove
  the convergence of the many-particle entropy towards the Boltzmann
  entropy of the solution to the limit equation (microscopic
  justification of the $H$-theorem of Boltzmann in this
  context). These results crucially rely on a new theory of
  quantitative uniform in time estimates of propagation of chaos. 

\vskip 0.5\baselineskip

\selectlanguage{francais}
\noindent{\bf R\'esum\'e} \vskip 0.5\baselineskip \noindent {\bf
  \`A propos du Programme de Kac en Th\'eorie Cin\'etique. }
Dans cette Note, nous pr\'esentons les r\'esultats principaux du
travail r\'ecent \cite{mm-article}, qui r\'epond \`a plusieurs
conjectures propos\'ees il y a une cinquantaine d'ann\'ees par Kac
\cite{kac}. Dans ce travail Kac introduit un processus stochastique
\`a grand nombre de particules (aujourd'hui appel\'e \'equation
ma\^itresse de Kac) qui converge, pour des donn\'ees chaotiques, vers
l'\'equation de Boltzmann spatialement homog\`ene. Nous r\'epondons
aux trois questions suivantes soulev\'ees dans cet article~: (1)
prouver la propagation du chaos pour des processus de collision
r\'ealistes (dans notre cas~: sph\`eres dures et \og vraies\fg{}
mol\'ecules maxwelliennes), (2) connecter les vitesses de relaxation
du processus stochastique et de l'\'equation limite en obtenant des
taux ind\'ependants du nombre de particules, (3) prouver la
convergence de l'entropie en grand nombre de particules vers
l'entropie de Boltzmann pour la solution de l'\'equation limite
(justification microscopique du th\'eor\`eme $H$ dans ce
contexte). Tous ces r\'esultats font appel de mani\`ere cruciale \`a
une nouvelle th\'eorie d'estimations quantitatives et uniformes en
temps de propagation du chaos.

\end{abstract}

\end{frontmatter}

\selectlanguage{francais}
\section*{Version fran\c{c}aise abr\'eg\'ee}

Le programme de Kac en th\'eorie cin\'etique consiste \`a comprendre
comment d\'eduire l'\'equation de Boltzmann spatialement homog\`ene
\`a partir d'un processus stochastique de saut sur l'espace des
vitesses \`a grand nombre de particules. Le but de ce programme est de
comprendre la notion de \og chaos mol\'eculaire\fg{} dans un cadre
plus simple que celui de la dynamique compl\`ete des particules, ainsi
que de donner une justification microscopique au th\'eor\`eme $H$
(croissance de l'entropie) et au processus de retour vers
l'\'equilibre. Nous renvoyons \`a la version compl\`ete pour
l'introduction de l'\'equation de Boltzmann \eqref{el}, des processus
de saut consid\'er\'es, ainsi que pour les d\'efinitions de la notion
de chaos et de la distance de Wasserstein $W_1$.  \medskip

\begin{theoreme}[R\'esum\'e des r\'esultats principaux]\label{theo:main:fr}
  On consid\`ere $d \ge 2$ et une distribution initiale $f_0 \in
  P(\R^d) \cap L^\infty$ \`a support compact ou poss\'edant
  suffisamment de moments polyn\^omiaux born\'es, et que l'on suppose
  centr\'ee sans perte de g\'en\'eralit\'e. Soit $f_t$ la solution
  correspondante de l'\'equation de Boltzmann \eqref{el} pour les
  sph\`eres dures ou les mol\'ecules maxwelliennes (sans troncature
  angulaire), et soit $f^N_t$ la solution du processus de saut \`a $N$
  particules correspondant, avec pour donn\'ee initiale $f_0 ^N$: soit
  (a) la tensorisation de $f_0 ^{\otimes N}$ de $f_0$, ou (b) la
  tensoris\'ee $f^{\otimes N} _0$ conditionn\'ee \`a la sph\`ere
  $\mathcal S^N$ (d\'efinie par \eqref{sphere-bol}).

  \begin{enumerate}
  \item {\bf Propagation de chaos quantifi\'ee et uniforme en temps}: 
On consid\`ere le cas (a) o\`u $f_0 ^N = f_0 ^{\otimes N}$. Alors
\begin{equation*}
  \forall \, N \ge 1, \ \forall \, 1 \le \ell  \le N, \quad 
  \sup_{t \ge 0} {W_1 \left( \Pi_{\ell} f^N
      _t, \left( f_t ^{\otimes \ell} \right) \right) \over \ell} \le \alpha(N)
\end{equation*}
avec $\alpha(N) \to 0$ lorsque $N \to \infty$, et o\`u $\Pi_\ell g^N$
d\'esigne la $\ell$-marginale d'une probabilit\'e $g$
sur $(\mathbb R^d)^N$.

\item {\bf Propagation du chaos entropique}: On consid\`ere le cas (b)
  o\`u $f_0 ^N$ est conditionn\'ee \`a $\mathcal S^N$. Alors 
la solution est entropie-chaotique~:
$$
\forall \, t \ge 0, \quad \frac1N \, H\left( f^N _t | \gamma^N \right)
\to H\left(f_t | \gamma \right), \quad N \to +\infty
$$
(voir \eqref{rel-entropie} pour les d\'efinitions des fonctionnelles
$H$) avec $\gamma$ la probabilit\'e gaussienne centr\'ee d'\'energie
$\mathcal E$ \'egale \`a l'\'energie de $f_0$ et $\gamma^N$ la mesure
de probabilit\'e uniforme sur $\mathcal S^N$. Cela fournit une
d\'erivation microscopique du th\'eor\`eme $H$ dans ce contexte.

\item {\bf Taux de relaxation ind\'ependants du nombre de particules}:
  On consid\`ere le cas (b) o\`u $f_0 ^N$ est conditionn\'ee \`a $\mathcal
  S^N$. Alors
\begin{equation*}
  \forall \, N \ge 1, \ \forall \, 1 \le \ell  \le N, \  \forall \, t
  \ge 0, \quad  
    {W_1 \left( \Pi_{\ell} f^N
      _t, \Pi_\ell \left( \gamma^N \right) \right) \over \ell} \le
  \beta(t) \qquad \mbox{avec} \ \beta(t) \to 0, \ t \to 0. 
\end{equation*}
Dans le cas des mol\'ecules maxwelliennes, et si la donn\'ee initiale
poss\`ede une information de Fisher finie (voir \eqref{fisher}), on
prouve \'egalement 
$
\forall \, N \ge 1, \ 0 \le \frac1N \, H\left( f^N _t | \gamma^N \right)
\le \beta(t)  \quad \mbox{avec} \ \beta(t) \to 0, \ t \to 0.$ 
  \end{enumerate}
\end{theoreme}

\selectlanguage{english}
\setcounter{equation}{0}
\section{Introduction}\label{Sec:Intro}

Motivated by the understanding of irreversibility and ``molecular
chaos'' in the Boltzmann equation
\cite{Maxwell1867,Boltzmann1872,Boltzmann1896} Mark Kac proposed in
1956 \cite{kac,kac2} the simpler and seemingly more tractable question
of deriving the \emph{spatially homogeneous} Boltzmann equation from a
\emph{many-particle jump proces}s, and he introduced a rigorous notion of
molecular chaos in this context. He proposed the first proof of the
\emph{propagation of chaos} along time for a simplified collision process for
which series expansions of the solution are available, and he showed
how the many-particle limit rigorously follows from it. 
A key motivation for Kac was the microscopic derivation of the
$H$-theorem (monotonicity of the Boltzmann entropy) in this context
which has remained open so far. Kac also raised the natural question
of connecting the asymptotic behavior of the many-particle process and
that of the limit nonlinear equation.  In his mind this program was to
be achieved by understanding dissipativity at the level of the linear
many-particle jump process and he insisted on the importance of
estimating its rate of relaxation. This has motivated beautiful works
on this ``Kac's spectral gap problem''
\cite{janvresse,maslen,CarlenCL2003}, but so far this strategy has
proved unsuccessful in obtaining relaxation rates which do not
degenerate in the many-particle limit, see the interesting discussion
in \cite{CCLLV}.

In this Note we present the main results in \cite{mm-article}. In this
paper we develop a quantitative theory of mean-field limit which
\emph{strongly relies on detailed knowledge of the limit nonlinear
  equation, rather than on detailed properties of the many-particle
  Markov process}. As the main outcome of this theory we prove uniform
in time quantitative propagation of chaos as well as propagation of
entropic chaos, and we prove relaxation rates \emph{independent of the
  number of particles} (measured in Wasserstein distance and relative
entropy). All this is done for the two important realistic and
achetypal models of collision, namely hard spheres and true (without
cutoff) Maxwell molecules. This provides a first complete answer to
the questions raised by Kac, however our answer is an ``inverse''
answer in the sense that our methodology is ``top-down'' from the
limit equation to the many-particle system rather than ``bottom-up''
as was proposed by Kac.

\subsection{The Boltzmann equation}
\label{sec:boltzmann-equation}


The {\it spatially homogeneous Boltzmann equation}  reads
 \begin{equation}\label{el}
   \frac{\partial f}{\partial t}(t,v) 
   = Q(f,f)(t,v), \qquad  v \in \R^d, \quad t \geq 0,
 \end{equation}
where $d \ge 2$ is the dimension and $Q$  is defined by
 \begin{equation*}
 Q(g,f)(v) = \frac12\,\int _{\R^d \times \mathbb{S}^{d-1}} B(|v-v_*|, \cos \theta)
         \left(g'_* f' + g' f_* '- g_* f - g f_* \right) \, dv_* \, d\sigma,
 \end{equation*}
where we have used the shorthands $f=f(v)$, $f'=f(v')$, $g_*=g(v_*)$ and
$g'_*=g(v'_*)$. Moreover, $v'$ and $v'_*$ are parametrized by
 \begin{equation*}
   v' = \frac{v+v_*}2 + \frac{|v-v_*|}2 \, \sigma, \qquad 
   v'_* = \frac{v+v_*}2 - \frac{|v-v_*|}2 \, \sigma, \qquad 
   \sigma \in \mathbb{S}^{d-1}. 
 \end{equation*}
 Finally, $\theta\in [0,\pi]$ is the deviation angle between $v'-v'_*$
 and $v-v_*$ defined by $\cos \theta = \sigma \cdot \hat u$, $u = v-v_*$, $\hat u = u/|u|$, 
  and $B$ is the Boltzmann collision
 kernel determined by physics (related to the cross-section
 $\Sigma(v-v_*,\sigma)$ by the formula $B=|v-v_*| \, \Sigma$).  


Boltzmann's collision operator has the fundamental properties of
conserving mass, momentum and energy
  \begin{equation*}
    \frac{d}{dt} \int_{\R^d} f \, \phi(v) \, dv = \int_{\R^d}Q(f,f) \, \phi(v)\,dv = 0, \quad
    \phi(v)=1,v,|v|^2, 
\end{equation*}
and satisfying the so-called Boltzmann's $H$ theorem which writes (at
the formal level)
\begin{equation*} 
  - \frac{d}{dt} H(f) := - \frac{d}{dt} \int_{\R^d} f \log f \, dv = -
  \frac{d}{dt} H(f |\gamma) := - \frac{d}{dt} \int_{\R^d} f \log
  \frac{f}{\gamma} \, dv = 
  - \int_{\R^d} Q(f,f)\log(f) \, dv \geq 0
\end{equation*}
where $\gamma$ is the gaussian with same mass, momentum and energy as
$f$. Note that the $H$ functional is the \emph{opposite} of the physical
entropy. 

We shall consider the folllowing important physical cases for $B$ (see~\cite{mm-article} for more details)
\[
B=\Gamma(|v-v_{*}|) \, b(\cos \theta) \qquad \mbox{ with } \ \Gamma,b \ge
0 \qquad \mbox{ given by one of the following formulas:} 
\]
 \begin{itemize}
 \item[(1)] {\bf (HS)} \textbf{Hard Spheres collision kernel}:
   $B(|v-v_*|, \cos \theta)= \Gamma(|v-v_*|) = 
    C \,
     |v-v_*|$ for some $C>0$.
\smallskip
\item[(2)] 
   {\bf (tMM)} \textbf{True Maxwell Molecules collision kernel}: \\ 
   $B(|v-v_*|, \cos \theta)= b(\cos \theta)
   \sim_{\theta \sim 0} C \, \theta^{-5/2}$ for some $C>0$.
 \smallskip

  \item[(3)] {\bf (GMM)} \textbf{Grad's cutoff Maxwell
     Molecules kernel}: 
   $B(|v-v_*|, \cos \theta)=1.$
 \end{itemize}
 

\subsection{Kac's program}
\label{sec:kacs-program}

{\em Kac's jump process} runs as follows: consider $N$ particles with
velocities $v_1$, \dots, $v_N \in \R^d$. Compute random times for each
pair of particles $(v_i,v_j)$ following an exponential law with
parameter $\Gamma(|v_i-v_j|)$, take the smallest, and perform a
collision $(v_i,v_j) \to (v_i ^*, v_j ^*)$ given by a random choice of
a direction parameter whose law is related to $b(\cos \theta)$, then
recommence. This process can be considered on $\R^{dN}$, however it
leaves invariant some submanifolds of $\R^{dN}$ (depending on the
number of conserved quantities during collision) and can be restricted
to them. In the original simplified model of Kac $d=1$ (scalar
velocities), the direction parameter is $\theta$ with collision rule
$$
v_i ^* = v_i \, \cos \theta + v_j \, \sin \theta, \qquad  
v_j ^* = - v_i \, \sin \theta + v_j \, \sin \theta
$$
and the collision process can be restricted to $\mathbb S^{N-1}(\sqrt{
  \mathcal E N})$ the sphere with radius $\sqrt{\mathcal E N}$, for
any given value of the \emph{energy} $\mathcal E$. For the more
realistic hard spheres of Maxwell molecules models, $d=3$, the
direction parameter is $\sigma \in \mathbb S^2$ with collision rule
$$ 
v_i ^* = \frac{v_i + v_j}2 + \frac{|v_i - v_j|}2 \, \sigma, \qquad  
v_j ^* = \frac{v_i + v_j}2 - \frac{|v_i - v_j|}2 \, \sigma \qquad 
\mbox{ with } \ \sigma \cdot \frac{(v_i-v_j)}{|v_i-v_j|} = \cos \theta 
$$
and the collision process can be restricted to the sphere
\begin{equation}\label{sphere-bol}
\mathcal S^N:= \mathbb S^{dN-1}\left( 
  \sqrt{N \mathcal E}\right) \cap \left\{ v_1 + \dots + v_N =0 \right\}.
\end{equation}

Kac formulated the notion of {\em propagation of chaos} that we shall
now explain. Consider a sequence $(f^N)_{N \ge 1}$ of probabilities on
$\mathbb R^{dN}$: the sequence is said $f$-\textit{chaotic} if $f^N
\sim f^{\otimes N}$ when $N \to \infty$ for some given one-particle
probability $f$ on $\mathbb R^d$. The meaning of this convergence is the
following: convergence in the weak measure topology for any marginal
depending on a finite number of variables. This is a \emph{low
  correlation} assumption.  It was clear since Boltzmann that in the
case when the joint probability density $f^N$ of the $N$-particle
system is tensorized \emph{during some time interval} into $N$ copies
$f^{\otimes N}$ of a $1$-particle probability density, then the latter
would satisfy the limit nonlinear Boltzmann equation during this time
interval. In general interactions between particles prevent any
possibility of propagation of the ``tensorization'' property, however
if the weaker property of chaoticity can be propagated along time in
the correct scaling limit it is sufficient for deriving the limit
equation.
Kac hence proved the propagation of chaos (with no rate) on the
simplified collision rule above (with $B=1$). His beautiful
combinatorial argument is based on an infinite series ``tree''
representation of the solution according to the collision history of
particles, and a Leibniz derivation-like formula for the iterated
$N$-particle operator acting on tensor products.

He then raises several questions that we schematize as follows:
\begin{enumerate}
\item The first one is concerned with the restriction of the models as
  compared to realistic collision processes: {\bf can one prove
    propagation of chaos for the hard spheres collision process?}

\item Following closely the spirit of the previous question 
  it seems to us
  very natural to ask whether {\bf one can prove propagation of chaos
    for the true Maxwell molecules collision process}? 
  This is related with \emph{long-range interactions} and
  \emph{fractional derivative operators}.

\item 
  Kac {\bf conjectures the propagation of the convergence of the
    $N$-particle $H$-functional towards the limit $H$-functional of
    the solution to the limit equation along time in the mean-field
    limit}. Since the latter always decays for a many-particle jump
  process, in his words ``{\it If the above steps could be made
    rigorous we would have a thoroughly satisfactory justification of
    Boltzmann's $H$-theorem.}'' 

\item He finally discusses the relaxation times, with the goal of
  deriving relaxation times of the limit equation from the
  many-particle system. This imposes to have estimates
  \emph{independent of the number of particles} on this relaxation
  times: 
  {\bf can one
    prove relaxation times {\it independent of the number of
      particles} in Wasserstein distance and/or 
    relative entropy?}
\end{enumerate}

This paper is concerned with solving the four questions
outlined above.

\section{Main results}\label{Sec:Main}

\subsection{A few words on previous results}
\label{sec:review}

For Boltzmann collision processes, Kac \cite{kac} has proved the
propagation of chaos in the case of his baby one-dimensional model. It
was generalized by McKean \cite{McKean1967} to the Boltzmann collision
operator for ``Maxwell molecules with cutoff'', i.e. the case {\bf
  (GMM)} above (see also \cite{T2} for a partial result for non-cutoff
Maxwell molecules). Gr\"unbaum~\cite{Grunbaum} then proposed in a very
compact and abstract paper another method for dealing with hard
spheres, based on the Trotter-Kato formula for semigroups and a clever
functional framework. Unfortunately this paper was incomplete for
several reasons (see the discussion in \cite{mm-article}). A
completely different approach was undertaken by Sznitman in the
eighties \cite{S1,S6} and he gave a full proof of propagation of chaos
for hard spheres by a probabilistic (non-constructive) approach.  Let
us also emphasize several quantitative results on a finite time
interval by Graham, M\'el\'eard and Fournier for Maxwell molecules
models \cite{GM,FM7,FM10}, and the works on the so-called \og Kac's
spectral gap problem\fg{} \cite{janvresse,maslen,CarlenCL2003,CCLLV}.

\subsection{Main results}
\label{sec:intromainresults}

\begin{theorem}[Summary of the main results]\label{theo:main}
  Consider some initial distribution $f_0 \in P(\R^d) \cap L^\infty$
  with compact support or polynomial moment bounds, taken to be
  centered without loss of generality. Consider the corresponding
  solution $f_t$ to the spatially homogeneous Boltzmann equation for
  hard spheres of Maxwell molecules, and the solution $f^N_t$ of the
  corresponding Kac's jump process starting either (a) from the
  tensorization $f_0 ^{\otimes N}$ of $f_0$ or (b) the latter
  conditionned to $\mathcal S^N$ (defined in (\ref{sphere-bol})).

The results in \cite{mm-article} can be classified into three main statements:
  \begin{enumerate}
  \item {\bf Quantitative uniform in time propagation of chaos (with
      any number of marginals)}:
\begin{equation*}
  \forall \, N \ge 1, \ \forall \, 1 \le \ell  \le N, \quad 
  \sup_{t \ge 0} {W_1 \left( \Pi_{\ell} f^N
      _t, \left( f_t ^{\otimes \ell} \right) \right) \over \ell} \le \alpha(N)
\end{equation*}
for some $\alpha(N) \to 0$ as $N \to \infty$, where $\Pi_\ell g^N$
stands for the $\ell$-marginal of an $N$-particle distribution $g^N$,
and where $W_1$ is the Wasserstein distance between probabilities on
$\mathbb R^{d\ell}$:
$$
W_1(p_1,p_2) := \sup_{[ \varphi ]_{\mbox{{\tiny {\em Lip}}}(\mathbb
    R^{d\ell})} \le 1} \int_{\mathbb R^{d\ell}} \varphi \, (dp_1 - dp_2)
\qquad \mbox{where } \ [\cdot]_{\mbox{{\tiny {\em Lip}}}} \ \mbox{
  denotes the Lipschitz semi-norm.}
$$ 
In the case (a) $f^N _0 = f^{\otimes N} _0$ one has moreover explicit
power law rate (for Maxwell molecules) or logarithmic rate (for hard
spheres) estimates on $\alpha$.

\item {\bf Propagation of entropic chaos}: Consider the case (b) where the
  initial datum of the many-particle system is restricted to $\mathcal S^N$. 
  Then if the initial datum is entropy-chaotic in the sense
$$
 \frac1N \, H\left( f^N _0 | \gamma^N \right) \to H\left(f_0 | \gamma
\right), \quad N \to +\infty
$$
\begin{equation}\label{rel-entropie}
\mbox{with} \quad H\left( f^N _0 | \gamma^N \right) := \int_{\mathcal
  S^N} \frac{df^N _0}{d \gamma^N} \, \log \frac{df^N
  _0}{d\gamma^N} \, \gamma^N(dV) \ \mbox{ and } \ H\left(f_0 | \gamma \right) :=
\int_{\R^d} f_0 \, \log \frac{f_0}{\gamma} \, dv
\end{equation}
and where $\gamma$ is the gaussian equilibrium with energy $\mathcal
E$ and $\gamma^N$ is the uniform probability measure on $\mathcal
S^N$, then the solution is also entropy-chaotic for any later time:
$$
\forall \, t \ge 0, \quad \frac1N \, H\left( f^N _t | \gamma^N \right)
\to H\left(f_t | \gamma \right), \quad N \to +\infty.
$$
Since our $f_0 ^N$ is entropy-chaotic, this proves the
derivation of the $H$-theorem in this context.

\item {\bf Quantitative estimates on relaxation times, independent of
    the number of particles}: Consider the case (b) where the
  initial datum of the many-particle system is restricted to $\mathcal
  S^N$. Then
\begin{equation*}
  \forall \, N \ge 1, \ \forall \, 1 \le \ell  \le N, \  \forall \, t
  \ge 0, \quad  
    {W_1 \left( \Pi_{\ell} f^N
      _t, \Pi_\ell \left( \gamma^N \right) \right) \over \ell} \le
  \beta(t)  \qquad \mbox{with} \ \beta(t) \to 0, \ t \to 0. 
\end{equation*}
Moreover in the case of Maxwell molecules, and assuming moreover that
the Fisher information of the initial datum $f_0$ is finite:
  \begin{equation}\label{fisher}
  \int_{\R^d} \frac{\left| \nabla_v f_0 \right|^2}{f_0} \, dv <
  +\infty,
\end{equation}
the following estimate also holds:
$
\forall \, N \ge 1, \quad 0 \le \frac1N \, H\left( f^N _t | \gamma^N \right)
\le \beta(t)  \quad \mbox{with} \ \beta(t) \to 0, \ t \to 0. 
$ 
 \end{enumerate}
\end{theorem}

\section{A few words on the methods and proofs}\label{Sec:Proofs}

Let us briefly explain some ideas underlying the result
Theorem~\ref{theo:main}-(i). The other results are then obtained on
the basis of this key estimate, combined with the other latest
results obtained in this field. 
\begin{itemize}
\item We aim at reducing the problem to a stability analysis when
  approximating a linear semigroup. To this purpose a key idea is to
  compare the linear $N$-particle dual evolution in $C_b(\mathbb
  R^{dN})$ with the (linear!) \emph{push-forward} evolution associated
  with the limit equation.
\item This push-forward semigroup is defined as follows: if $S^{N
    \!L}_t$ denotes the nonlinear semigroup of \eqref{el}, this
  push-forward semigroup is defined on $P(P(\mathbb R^d))$ by
  $T^\infty _t[\Phi](f) = \Phi(S^{N \!L}_t(f))$.
\item In order to make this comparison between semigroups, we use the
  empirical measure $\mu^N _V = (\sum_{i=1} ^N \delta_{v_i})/N$ in
  order to embed the dynamics in $P(\mathbb R^{dN})$ into a dynamics
  in $P(P(\mathbb R^d))$.
\item One then considers the following term to be estimated 
$$
\left| \left \langle \left( S^N_t(f_0^{N}) - \left( S _t ^{NL}
          (f_0)
        \right)^{\otimes N} \right),  \varphi  
      \otimes 1^{\otimes N-\ell} \right\rangle \right|
$$
for some test function $\varphi$ only depending on $\ell$ variables. 

\item The approximation of these marginals by empirical
  measure estimates yields a first error term on the $N$-particle
  semigroup
$$
\left| \left\langle S^N_t(f_0^N), 
      \varphi  \otimes
      1^{\otimes N-\ell} \right\rangle -
    \left \langle S^N_t(f_0^N), 
      R^\ell_\varphi \circ \mu^N_V \right\rangle  \right| \qquad
  \mbox{ with } \ \ R^\ell _\varphi (f) := \int_{\mathcal R^{d\ell}}
  \varphi \, f^{\otimes \ell}(dv_1 \dots dv_\ell)
$$
which is controlled by combinatorial arguments, and then a second
error term 
$$
\left| \left\langle f_0^N, 
      (T_t ^\infty R^\ell_\varphi ) \circ \mu^N_V) \right\rangle
    -  \left\langle (S _t ^{NL} (f_0))^{\otimes \ell} ,
      \varphi \right\rangle \right|
$$
which is controlled thanks a stability for measure solutions of the
limit equation. 

\item Finally there remains the most important term where the two
  dynamics are effectively compared 
$$
\left| \left\langle f_0^N, T^N_t ( R^\ell_\varphi \circ
      \mu^N_V) \right\rangle
    - \left\langle f_0^N, 
      (T_t ^\infty R^\ell_\varphi ) \circ \mu^N_V) \right\rangle
  \right|.
$$
This term is controlled by using (1) a quantitative argument \`a la
Trotter-Kato in order to express the difference of semigroups in terms
of the difference of their generators $G^N$ and $G^\infty$, (2) a
\emph{consistency estimate} between those generators, (3) a \emph{stability
estimate} on the limit equation. 
\item There is a \emph{loss of derivative} in the consistency estimate
  in the sense of differentiable functions acting on $P(\mathbb R^d)$,
  which lead us to develop a differential calculus on this space
  adapted to our purpose.
\item The stability estimate means, once translated on the original
  nonlinear semigroup $S^{N \!L}_t$ of \eqref{el}, a propagation of a
  bound $C^{1+\theta}(P(\mathbb R^d))$ on $S^{N \!L}_t$. The role
  played by such stability estimates is a key novelty of our
  study. Proving them for Boltzmann is also one of the most technical
  aspects of \cite{mm-article}. 
\item There are many possible choices of distances on the space of
  probabilities (total variation but also many non-equivalent weak
  measure distances), and it is a crucial point that our method is
  flexible enough to allow for many such different choices adapted to
  the equations it is applied to. 

\end{itemize}



\end{document}